\documentclass[11pt,a4paper]{amsart}
\usepackage[T1]{fontenc}
\usepackage{lmodern}
\usepackage[margin=26mm]{geometry}
\usepackage{amsmath,amssymb,mathtools}
\usepackage[all]{xy}
\usepackage{microtype}
\usepackage{needspace}
\usepackage[hidelinks]{hyperref}
\hypersetup{pdftitle={Cartesian closedness of lax comma categories},pdfauthor={Maria Manuel Clementino and Fernando Lucatelli Nunes}}
\newtheorem{theorem}{Theorem}[section]
\newtheorem{lemma}[theorem]{Lemma}
\newtheorem{proposition}[theorem]{Proposition}
\newtheorem{corollary}[theorem]{Corollary}
\theoremstyle{definition}
\newtheorem{definition}[theorem]{Definition}
\newtheorem{remark}[theorem]{Remark}

\numberwithin{equation}{section}
\newcommand{\Cat}{\mathsf{Cat}}
\newcommand{\CAT}{\mathsf{CAT}}
\newcommand{\Set}{\mathsf{Set}}
\newcommand{\Fam}{\mathsf{Fam}}
\newcommand{\CMon}{\mathsf{CMon}}
\newcommand{\Vect}{\mathsf{Vect}}
\newcommand{\Ab}{\mathsf{Ab}}
\newcommand{\Gr}{\mathsf{Gr}}
\newcommand{\Disc}{\mathsf{Disc}}
\newcommand{\Ob}{\mathsf{Ob}}
\newcommand{\one}{\mathbf{1}}
\newcommand{\id}{\mathrm{id}}
\newcommand{\op}{\mathrm{op}}
\newcommand{\sslash}{\mathbin{/\mkern-6mu/}}
\newcommand{\DX}{\Cat\sslash X}
\newcommand{\DY}{\Cat\sslash Y}

\newcommand{\Adj}[4]{\xymatrix@C=4.5em{#1\ar@/_1.5pc/[r]_{#4}& #2\ar@/_1.5pc/[l]_{#3}\ar@{}[l]|{\bot}}}
\title[Cartesian closedness of lax comma categories]{Cartesian closedness of lax comma categories}
\author{Maria Manuel Clementino}
\address{CMUC, Department of Mathematics, University of Coimbra, 3000-143 Coimbra, Portugal}
\email{mmc@mat.uc.pt}
\author{Fernando Lucatelli Nunes}
\address{CMUC, Department of Mathematics, University of Coimbra, 3000-143 Coimbra, Portugal}
\email{fln@uc.pt}
\date{}
\begin{document}
\begin{abstract}
For a finitely complete category $X$, we prove that the lax comma category $\Cat\sslash X$ is cartesian closed if and only if its free coproduct completion $\Fam(X)$ is cartesian closed. The proof uses a change-of-base theorem: postcomposition with a functor induces a left adjoint between lax comma categories precisely when that functor is a left multiadjoint. We describe the exponentials and derive examples from known results on family categories.
\end{abstract}
\maketitle

\section*{Introduction}
The study of lax comma categories in~\cite{CN2} was motivated by
two-dimensional aspects of Janelidze-Galois theory~\cite{BJ} and
their connection with lax orthogonal factorization systems.
Subsequent work on lax comma categories of ordered sets~\cite{CNord}
and categories~\cite{CNP} addressed exponentiability, extensivity,
topologicity and descent. Related developments include the study of
filtered preorders and ordered families~\cite{CJ,CP}, as well as
topological lax comma categories~\cite{CHP}. The latter also occur
in the description of the arithmetic spaces used in the Stone--Gelfand
duality of Abbadini, Marra and Spada~\cite{AMS}.

The cartesian closedness result of~\cite[Theorem~2.4]{CNP} constructs exponentials in $\Cat\sslash X$ using exponentials and ends in the base. In particular, a complete cartesian closed $X$ gives a cartesian closed $\Cat\sslash X$, and the projection to $\Cat$ preserves these exponentials. However, the results on $\Sigma$-tractable monoidal structures in~\cite{NV} show that $\Cat\sslash\CMon$ is cartesian closed, although $\CMon$ is not. Thus the fibred construction does not exhaust the possibilities. The challenge is to characterize cartesian closedness of the total category $\Cat\sslash X$ without requiring its exponentials to be preserved by the projection
$ \Cat\sslash X\to \Cat .$\\
We answer this question in terms of families, that is to say, the free coproduct completion of $X$. For finitely complete $X$, we prove
\[
 \DX\text{ is cartesian closed}
 \quad\Longleftrightarrow\quad
 \Fam(X)\text{ is cartesian closed}.
\]
The proof rests on a general change-of-base result: left multiadjoints are precisely the functors inducing left adjoints between lax comma categories. For the singleton embedding, this gives a full coreflection of $\DX$ in $\Cat\sslash\Fam(X)$. Together with the discrete-diagram embedding, it allows cartesian closedness to be transferred between families (objects of $\Fam(X)$) and diagrams (objects of $\DX$).

Cartesian closedness of $\Fam(X)$ has already been studied in several settings, from the characterization of Ad\'amek and Rosick\'y~\cite{AR} to the semantics of automatic differentiation~\cite{CHAD} and free distributive and extensive completions~\cite{FED,Dist}. The final section applies these results to lax comma categories.

\section{Families}\label{sec:families}

Throughout, families and diagram shapes are small, while the categories $X$ and $Y$ may be large. Presheaf categories and the ordinary comma constructions recalled in Section~\ref{sec:diagrams} may be interpreted in a larger universe. We identify sets with discrete categories and denote the singleton set by $\one$. An object $x\in X$ also denotes the functor $\one\to X$ selecting it. When $X$ has a terminal object, we denote it by $1$.

For limits and cartesian closedness in the free coproduct completion $\Fam(X)$, see Hu and Tholen~\cite{HT} and Ad\'amek and Rosick\'y~\cite{AR}, as well as~\cite[Sections~6.1 and~6.2]{BJ}. Family categories occur in the semantics of automatic differentiation~\cite{CHAD,ICHAD}, monoidal closure of Grothendieck constructions~\cite{NV}, and free distributive and extensive completions~\cite{FED,Dist}. Descent in categories of families is studied in~\cite{Prezado,CNP,CP}.

An object of the free coproduct completion $\Fam(X)$ is a pair $(A,F)$, where $A$ is a set and $F:A\to X$ is a functor. A morphism $(u,\alpha):(A,F)\to(B,G)$ consists of a function $u:A\to B$ and a natural transformation $\alpha:F\Rightarrow G\circ u$; thus its components are morphisms $\alpha_a:F(a)\to G(u(a))$. The singleton embedding $j_X:X\to\Fam(X)$, given by $j_X(x)=(\one,x)$, is fully faithful. We often write $j$ for $j_X$.

Coproducts are given by disjoint unions of the indexing sets. If $X$ has finite products, then $j(1)=(\one,1)$ is terminal in $\Fam(X)$ and
\begin{equation}\label{eq:products}
 (A,F)\times(B,G)=(A\times B,F\times G),
 \qquad (F\times G)(a,b)=F(a)\times G(b).
\end{equation}
In particular, $j$ preserves finite products. More generally, if the required products exist in $X$, then
\[
 \prod_{k\in K}(A_k,F_k)=\left(\prod_{k\in K}A_k,H\right),
 \qquad H(s)=\prod_{k\in K}F_k(s(k)).
\]
If $X$ has finite limits, so does $\Fam(X)$; see~\cite[Remarks~2.2(1) and~2.3]{AR}.

The functor
\begin{equation}\label{eq:yoneda}
 \Phi_X:\Fam(X)\longrightarrow\CAT[X^{\op},\Set],
 \qquad (A,F)\longmapsto\coprod_{a\in A}X(-,F(a))
\end{equation}
is fully faithful by Yoneda, with essential image the \emph{familially representable} presheaves~\cite{CJohnstone}. We recall Diers's multi-notions~\cite{Diers} in this language; see also~\cite{AR,HT}.

\begin{definition}\label{def:multi}
A \emph{multilimit} of $M:A\to X$ is an object $(I,L)\in\Fam(X)$ representing its cone presheaf:
\[
 \Fam(X)(j(z),(I,L))\cong\CAT[A,X](\Delta(z),M).
\]
The category $X$ is \emph{multicomplete} if every diagram has a multilimit. For discrete diagrams we speak of \emph{multiproducts}. If $X$ has finite products, it is \emph{cartesian multi-closed} if each functor $X(-\times x,y)$ is familially representable.

A functor $c:X\to Y$ is a \emph{left multiadjoint} if, for every $y\in Y$, the presheaf $Y(c(-),y)$ is familially representable. Thus there is a family $\widehat R(y)\in\Fam(X)$ and an isomorphism
\begin{equation}\label{eq:multiadjoint}
 Y(c(x),y)\cong\Fam(X)(j(x),\widehat R(y)),
\end{equation}
natural in $x$. Equivalently, $(c\downarrow y)$ has a set of connected components, each with a terminal object.
\end{definition}

By full faithfulness of~\eqref{eq:yoneda}, the families $\widehat R(y)$ extend to a functor $\widehat R:Y\to\Fam(X)$, with~\eqref{eq:multiadjoint} natural in both variables. Equivalently, $\Fam(c)$ has a right adjoint: its value at $(B,G)$ is $\coprod_{b\in B}\widehat R(G(b))$. This is the family-completion formulation of multiadjointness; see~\cite[Section~2.6]{AST} and~\cite[Remark~2.7]{AR}.

\begin{lemma}[{\cite[Remark~2.2(1) and Lemma~2.6]{AR}}]\label{lem:multicomplete}
The category $\Fam(X)$ is complete if and only if $X$ is multicomplete. It has all products if and only if $X$ has all multiproducts.
\end{lemma}

\section{Lax comma categories and change of base}\label{sec:diagrams}

We use the convention of~\cite{CNP}. An object of $\DX$ is a pair $(A,F)$, where $A$ is a small category and $F:A\to X$. A morphism $(A,F)\to(B,G)$ is a pair $(u,\alpha)$ with
\[
 u:A\longrightarrow B,
 \qquad \alpha:F\Longrightarrow G\circ u.
\]
Composition is $(v,\beta)\circ(u,\alpha)=(v\circ u,\gamma)$, where $\gamma_a=\beta_{u(a)}\circ\alpha_a$. When $X$ has finite products, $(\one,1)$ is terminal in $\DX$ and
\begin{equation}\label{eq:diagproducts}
 (A,F)\times(B,G)=(A\times B,F\times G),
 \qquad (F\times G)(a,b)=F(a)\times G(b).
\end{equation}

For the recalled adjunction~\eqref{eq:strictadj}, Remark~\ref{rem:diagrams}, and the comparison~\eqref{eq:comparison}, we work in an ambient universe containing the comma objects and enlarge the slice and diagram categories accordingly.

For $c:X\to Y$, write $c_!:\DX\to\DY$ for postcomposition and $I_X:\Cat/X\to\DX$ for the inclusion of the strict slice. The comma construction of~\cite[Theorem~5.8]{CN2}, in its $2$-cell-dual form, gives the adjunction
\begin{equation}\label{eq:strictadj}
 \Adj{\DY}{\Cat/X}{c_!\circ I_X}{c^{\Rightarrow}},
 \qquad c^{\Rightarrow}(B,G)=(c\downarrow G\longrightarrow X).
\end{equation}
Indeed, a strict map into the comma object consists of a functor $u:A\to B$ and a transformation $c\circ F\Rightarrow G\circ u$. In particular, $I_X\dashv K_X$, where
$K_X(B,G)=(\id_X\downarrow G\to X)$.

\begin{remark}\label{rem:diagrams}
The notation $\Cat\sslash X$ and the above definition, with the same direction of the comparison transformation, already appear in Grothendieck~\cite[Section~11(b), p.~186]{SGA1}. This is also the category of small diagrams considered by Peschke and Tholen~\cite[Definition~2.1]{PeschkeTholen} and Perrone and Tholen~\cite[Definition~2.1]{PerroneTholen}. Its projection to $\Cat$ is the split fibration associated with $A\mapsto\CAT[A,X]$, as recalled in~\cite[Section~1.2]{CNP}.

In the ambient universe, the adjunction $I_X\dashv K_X$ induces the free split fibration $2$-monad on the strict slice:
\[
 T_X(F)=\bigl(\id_X\downarrow F\xrightarrow{\pi_X}X\bigr).
\]
Its unit inserts identity arrows and its multiplication composes arrows in $X$. Its strict algebras are split fibrations over $X$; see Street~\cite[Section~2]{Street74}. The corresponding comma $2$-monad with varying base and its strict algebras are considered in~\cite{NT}. A strict morphism $F\to T_X(G)$ over $X$ is precisely a pair $(u,\alpha)$ with $\alpha:F\Rightarrow G\circ u$. Consequently, the underlying Kleisli category, restricted to objects with small domains, is $\DX$. This is the $2$-cell-dual form of~\cite[Theorem~5.11]{CN2} and gives the connection with free fibrations recorded in~\cite[Section~4.3]{CNP}.
\end{remark}

We ask when $c_!$ itself has a right adjoint. Its construction uses a diagram whose domain records the indices in a family-valued functor.

For $T:B\to\Fam(X)$, write $T(b)=(S(b),T_b)$, with $T_b:S(b)\to X$, and $T(v)=(S(v),\tau_v)$. Thus $S:B\to\Set$ is a functor and $\tau_v:T_b\Rightarrow T_{b'}\circ S(v)$ for $v:b\to b'$. Put
\begin{equation}\label{eq:flatten}
 \Gr(B,T)=\left(\int_B S,\widetilde T\right),
 \qquad \widetilde T(b,s)=T_b(s).
\end{equation}
An arrow $(b,s)\to(b',s')$ is an arrow $v:b\to b'$ with $S(v)(s)=s'$; its image is $(\tau_v)_s:T_b(s)\to T_{b'}(s')$. The category $\int_B S$ is the Grothendieck construction of $S$, regarded as a functor with discrete values, and its projection to $B$ is a discrete opfibration. For the Grothendieck construction and the corresponding fibred and cofibred descriptions, see Grothendieck~\cite[Section~8]{SGA1}, Gray~\cite{Gray} and Street~\cite[Section~1.10]{Street}; a recent account is given in~\cite[Section~2.1]{MV}.

On a lax morphism $(v,\gamma):(B,T)\to(B',T')$, the index maps of $\gamma$ induce the functor between categories of elements, and its component maps give the comparison transformation. Thus $\Gr$ is functorial.

\begin{theorem}\label{thm:multiadjoint}
For a functor $c:X\to Y$, the following conditions are equivalent:
\begin{enumerate}
\item $c$ is a left multiadjoint;
\item $\Fam(c)$ has a right adjoint;
\item $c_!:\DX\to\DY$ has a right adjoint.
\end{enumerate}
If $\widehat R$ is as in~\eqref{eq:multiadjoint}, the right adjoint in \textup{(3)} is
\[
 R(B,G)=\Gr(B,\widehat R\circ G).
\]
\end{theorem}
\begin{proof}
The equivalence of (1) and (2) was recalled in Section~\ref{sec:families}. Assuming (1), naturality of~\eqref{eq:multiadjoint} gives
\begin{align*}
 \DY((A,c\circ F),(B,G))
 &\cong\coprod_{u:A\to B}\CAT[A,Y](c\circ F,G\circ u)\\
 &\cong\coprod_{u:A\to B}\CAT[A,\Fam(X)](j\circ F,\widehat R\circ G\circ u)\\
 &\cong\DX((A,F),\Gr(B,\widehat R\circ G)).
\end{align*}
For the last bijection, write $\widehat R(G(b))=(S(b),E_b)$. A transformation $j\circ F\Rightarrow\widehat R\circ G\circ u$ selects $s_a\in S(u(a))$ and maps $F(a)\to E_{u(a)}(s_a)$. Naturality says precisely that $a\mapsto(u(a),s_a)$ is a functor to $\int_B S$ and that the component maps form a natural transformation. The bijections are natural, proving (3).

Conversely, suppose $c_!\dashv R$ and write $R(\one,y)=(D,E)$. Singleton diagrams give
\[
 Y(c(x),y)\cong\DX((\one,x),(D,E))
 \cong\coprod_{d\in\Ob(D)}X(x,E(d)),
\]
naturally in $x$. Since $D$ is small, this proves (1).
\end{proof}

We can now compare this adjunction with~\eqref{eq:strictadj}. Since $I_X\dashv K_X$, the composite $c_!\circ I_X$ has right adjoint $K_X\circ R$. Uniqueness of right adjoints gives, whenever the comma objects are available,
\begin{equation}\label{eq:comparison}
 c^{\Rightarrow}\cong K_X\circ R.
\end{equation}
Thus the adjunction of~\cite{CN2} factors through the two lax comma categories precisely when $c$ is a left multiadjoint.

The singleton embedding gives the following particular case.

\begin{corollary}\label{cor:flatten}
The functor $j_X:X\to\Fam(X)$ is a left multiadjoint. It induces a fully faithful left adjoint, preserving finite products whenever these exist in $X$:
\begin{equation}\label{eq:Jadj}
 \Adj{\Cat\sslash\Fam(X)}{\DX}{J=(j_X)_!}{\Gr}.
\end{equation}
\end{corollary}
\begin{proof}
For $(A,F)\in\Fam(X)$, the identity
$\Fam(X)(j(x),(A,F))=\coprod_{a\in A}X(x,F(a))$
shows that $j$ is a left multiadjoint with $\widehat R=\id_{\Fam(X)}$. By Theorem~\ref{thm:multiadjoint}, we have $J\dashv\Gr$. The unit is an isomorphism, since the index functor of $j\circ F$ is constantly a singleton. Product preservation follows from~\eqref{eq:products} and~\eqref{eq:diagproducts}.
\end{proof}

\begin{remark}
These adjunctions extend to $2$-adjunctions, in the usual $\Cat$-enriched sense~\cite{KS}. A $2$-cell $(u,\alpha)\Rightarrow(v,\beta)$ is a transformation $\tau:u\Rightarrow v$ such that $G(\tau_a)\circ\alpha_a=\beta_a$ for every $a\in A$. The bijections above identify precisely these transformations. The ordinary adjunctions suffice for cartesian closedness.
\end{remark}

\section{Cartesian closedness}

Families are diagrams with discrete domains: in this sense $\Fam(X)=\Set\sslash X$. The adjunction $\Disc\dashv\Ob$ between sets and categories induces
\begin{equation}\label{eq:discreteadj}
 \Adj{\DX}{\Fam(X)}{D}{U},
 \qquad U(A,F)=(\Ob(A),F\circ\iota_A),
\end{equation}
where $\iota_A:\Disc(\Ob(A))\to A$ is the canonical functor. Indeed, a morphism from a discrete diagram consists of an assignment of objects and a family of component maps. Thus $D$ is fully faithful and, when $X$ has finite products, preserves them. Together with~\eqref{eq:Jadj}, this gives the two full coreflections needed below.

\begin{lemma}\label{lem:coreflection}
Let $L:\mathcal A\to\mathcal B$ be a fully faithful left adjoint preserving finite products, with right adjoint $R$. If $L(P)$ is exponentiable in $\mathcal B$, then $P$ is exponentiable in $\mathcal A$, with $Q^P=R(L(Q)^{L(P)})$.
\end{lemma}
\begin{proof}
For every $Z\in\mathcal A$,
\[
 \mathcal A(Z,R(L(Q)^{L(P)}))
 \cong\mathcal B(L(Z)\times L(P),L(Q))
 \cong\mathcal A(Z\times P,Q).
\]
\end{proof}

\begin{proposition}\label{prop:necessary}
If $X$ has finite products and $\DX$ is cartesian closed, then $\Fam(X)$ is cartesian closed.
\end{proposition}
\begin{proof}
The functor $D$ in~\eqref{eq:discreteadj} is a fully faithful left adjoint preserving finite products. The conclusion therefore follows from Lemma~\ref{lem:coreflection}.
\end{proof}

For the converse, we use the products supplied by cartesian closedness of the family category.

\begin{theorem}\label{thm:lex}
Let $X$ be a category with finite limits. Then $\DX$ is cartesian closed if and only if $\Fam(X)$ is cartesian closed. Under these equivalent conditions, $\Fam(X)$ is complete.
\end{theorem}
\begin{proof}
We first observe that if $\Fam(X)$ is cartesian closed, then it is complete. Indeed, by~\cite[Theorem~2.11]{AR}, $X$ has all multiproducts, so $\Fam(X)$ has all products by Lemma~\ref{lem:multicomplete}. Since $X$ has finite limits, so does $\Fam(X)$. Products and equalizers therefore give all limits in $\Fam(X)$.

If $\DX$ is cartesian closed, Proposition~\ref{prop:necessary} gives cartesian closedness of $\Fam(X)$, and the preceding observation gives its completeness. Conversely, if $\Fam(X)$ is cartesian closed, it is also complete, so $\Cat\sslash\Fam(X)$ is cartesian closed by~\cite[Theorem~2.4]{CNP}. Corollary~\ref{cor:flatten} and Lemma~\ref{lem:coreflection} then give cartesian closedness of $\DX$.
\end{proof}

\section{Semiadditive case}\label{sec:further}

The known cartesian closedness of $\Fam(X)$ for complete semiadditive $X$ was one of our main motivations for investigating the non-fibred case for $\DX$. In fact, semiadditivity and products suffice~\cite[Proposition~74]{CHAD}; see also~\cite[Section~9]{NV}. The exponential formula in Corollary~\ref{cor:semiadditive} recovers the lax comma construction of~\cite[Section~9]{NV}.

The proof of Theorem~\ref{thm:lex} also describes the exponentials. Assume that $X$ has finite products and that $\Fam(X)$ is complete and cartesian closed. For $(A,F),(B,G)\in\DX$, define
\begin{equation}\label{eq:ends}
 T:\Cat[A,B]\longrightarrow\Fam(X),
 \qquad T(H)=\int_{a\in A}j(G(H(a)))^{j(F(a))}.
\end{equation}
The end is that of the bifunctor $(a,a')\mapsto j(G(H(a')))^{j(F(a))}$. Postcomposition makes it functorial in $H$. The construction of~\cite{CNP} and Lemma~\ref{lem:coreflection} give
\begin{equation}\label{eq:exponential}
 (B,G)^{(A,F)}=\Gr(\Cat[A,B],T).
\end{equation}
Writing $T(H)=(S(H),T_H)$, the domain of~\eqref{eq:exponential} is a discrete opfibration over $\Cat[A,B]$ with fibre $S(H)$ at $H$. This projection need not be an isomorphism.

The description is particularly simple for semiadditive categories. Recall that these are categories with finite biproducts. In this case
\begin{equation}\label{eq:biproduct}
 X(z\times x,y)\cong X(z,y)\times X(x,y)
 \cong\coprod_{\beta:x\to y}X(z,y).
\end{equation}
Thus $X$ is cartesian multi-closed~\cite[Example~2.9(b)]{AR}. If $X$ is complete,~\cite[Theorem~2.11]{AR} and Lemma~\ref{lem:multicomplete} also give the completeness and cartesian closedness of $\Fam(X)$.

\begin{corollary}\label{cor:semiadditive}
Let $X$ be a complete semiadditive category. Then $\DX$ is cartesian closed, with
\begin{equation}\label{eq:semiadditive}
 (B,G)^{(A,F)}=
 \left((F\downarrow G_!),\ (H,\beta)\longmapsto\lim_A(G\circ H)\right),
\end{equation}
where $G_!:\Cat[A,B]\to\CAT[A,X]$ is postcomposition with $G$.
\end{corollary}
\begin{proof}
The cartesian closedness assertion is an immediate corollary of Theorem~\ref{thm:lex} and the known cartesian closedness of $\Fam(X)$ for semiadditive categories with products~\cite[Proposition~74]{CHAD}; see also~\cite[Section~9]{NV}. To identify the exponential, for each $H:A\to B$ the biproduct decomposition gives
\[
 \CAT[A,X](\Delta(z)\times F,G\circ H)
 \cong\coprod_{\beta:F\Rightarrow G\circ H}X\left(z,\lim_A(G\circ H)\right).
\]
Thus $T(H)$ is the pair $(\CAT[A,X](F,G\circ H),\Delta\bigl(\lim_A(G\circ H)\bigr))$, where the second component is the constant functor. On a transformation $\eta:H\Rightarrow H'$, its index map sends $\beta$ to the transformation with components $G(\eta_a)\circ\beta_a$, and its component map is induced between limits. Flattening therefore gives exactly the comma category and diagram in~\eqref{eq:semiadditive}. Evaluation at $((H,\beta),a)$ is the biproduct map
$[\pi_a,\beta_a]:(\lim_A(G\circ H))\oplus F(a)\to G(H(a))$.
\end{proof}

For $X=\CMon$, the self-exponential of $(\one,\mathbb N)$ is $(\Disc(\mathbb N),\Delta(\mathbb N))$, with evaluation over $k$ given by $(m,n)\mapsto m+kn$. This recovers the non-fibred phenomenon in~\cite{NV}. The same construction applies to complete additive categories, including $\Ab$ and $\Vect_k$.

\section{Cartesian closed families}\label{sec:closedfamilies}

We conclude by relating our characterization to known results on family categories, starting with the following theorem of Ad\'amek and Rosick\'y.

\begin{theorem}[{\cite[Theorem~2.11]{AR}}]\label{thm:AR}
Let $X$ be a category with finite products. Then $\Fam(X)$ is cartesian closed if and only if $X$ has all multiproducts and is cartesian multi-closed.
\end{theorem}

Explicitly, let $R(x,y)\in\Fam(X)$ represent $X(-\times x,y)$. For $P=(A,F)$ and $Q=(B,G)$, one has
\begin{equation}\label{eq:famexp}
 Q^P=\prod_{a\in A}\coprod_{b\in B}R(F(a),G(b)).
\end{equation}
Indeed, a map from $j(z)$ into~\eqref{eq:famexp} selects a function $u:A\to B$ and maps $z\times F(a)\to G(u(a))$. This is a map $j(z)\times P\to Q$, and the assertion extends to all families by coproducts.

Together with Lemma~\ref{lem:multicomplete} and Theorem~\ref{thm:lex}, this gives the following formulation in terms of $X$.

\begin{corollary}\label{cor:intrinsic}
Let $X$ be a category with finite limits. Then $\DX$ is cartesian closed if and only if $X$ is multicomplete and cartesian multi-closed.
\end{corollary}

The examples of cartesian multi-closed categories in~\cite[Examples~2.9]{AR} include cartesian closed categories, semiadditive categories, and opposites of extensive categories. Whenever these are complete, Corollary~\ref{cor:intrinsic} applies. In particular, it gives cartesian closedness of $\Cat\sslash\Set^{\op}$, as well as the semiadditive examples above. Families of vector spaces and of opposite vector spaces also appear as cartesian closed semantic categories in \emph{CHAD for expressive total languages}~\cite[Section~10.2]{CHAD}; the present result gives their corresponding lax comma examples.

Free completions provide another source. Let $\mathcal L(C)$ denote the free completion of $C$ under limits. The cartesian closedness of $\Fam(\mathcal L(C))$ proved in~\cite[Theorem~4.4]{FED}, together with Theorem~\ref{thm:lex}, gives
\[
 \Cat\sslash\mathcal L(C)\quad\text{cartesian closed}.
\]
The same work treats exponentiable finite coproducts of connected objects in free infinitary lextensive categories~\cite[Theorem~4.7]{FED}.

Finally, the free doubly-infinitary distributive category $\Fam(\Fam(C^{\op})^{\op})$ is cartesian closed \cite[Theorem~3.11]{Dist}. Whenever $\Fam(C^{\op})^{\op}$ has finite limits, Theorem~\ref{thm:lex} therefore gives cartesian closedness of $\Cat\sslash\Fam(C^{\op})^{\op}$.

\subsection*{Acknowledgements}
The authors acknowledge partial financial support by the Centre for Mathematics of the University of Coimbra (CMUC), through the Funda\c{c}\~ao para a Ci\^encia e a Tecnologia (FCT), grants \href{https://doi.org/10.54499/UID/00324/2025}{UID/00324/2025} and UID/PRR/00324/2025.

\end{document}